\documentclass[11pt]{amsart}
\usepackage{amsthm,amsfonts,amssymb,euscript,setspace,enumitem,leftidx}

\newtheorem{thm}{Theorem}[section]
\newtheorem{cor}[thm]{Corollary}
\newtheorem{lem}[thm]{Lemma}

\theoremstyle{definition}

\theoremstyle{remark}

\numberwithin{equation}{section}

\newcommand{\R}{\mathbb{R}}
\newcommand{\Z}{\mathbb{Z}}

\newcommand{\F}{\mathcal{F}}

\newcommand{\X}{\mathcal{X}}

\newcommand{\ep}{\epsilon}
\newcommand{\la}{\lambda}
\newcommand{\La}{\Lambda}
\newcommand{\te}{\theta}
\newcommand{\al}{\alpha}
\newcommand{\be}{\beta}
\newcommand{\ga}{\gamma}

\newcommand{\de}{\delta}

\newcommand{\si}{\sigma}

\newcommand{\Om}{\Omega}

\def\beq{\begin{equation}}
\def\eeq{\end{equation}}

\newcommand{\p}{\partial}
\newcommand{\pd}[3][]{\frac{\p^{#1} #2}{\p^{#1} #3}}
\newcommand{\pdd}[3]{\frac{\p^2 #1}{\p #2\p #3}}
\newcommand{\lap}{\Delta}
\newcommand{\st}{\mid}

\newcommand{\closure}[1]{\overline{#1}}

\newcommand{\norm}[1]{\Vert #1 \Vert}
\newcommand{\abs}[1]{\left\vert #1 \right\vert}
\newcommand{\set}[1]{\left\{ #1 \right\}}

\newcommand{\hm}[1]{\mathcal{H}^{#1}} 
\newcommand{\compl}[1]{\leftidx{^c}{#1}{}}

\DeclareMathOperator{\dist}{dist}

\begin{document}

\title[Regularity of the Composite Membrane Problem]{Regularity of the Minimizers in the Composite Membrane Problem in $\R^2$}
\author{Sagun Chanillo}
\address{Department of Mathematics\\
Rutgers University\\
Piscataway, NJ 08854, USA}
\email{chanillo@math.rutgers.edu}

\author{Carlos E. Kenig}
\address{Department of Mathematics\\
University of Chicago\\
Chicago, IL 60637, USA}
\email{cek@math.uchicago.edu}

\author{Tung To}
\address{Department of Mathematics\\
University of Chicago\\
Chicago, IL 60637, USA}
\email{totung@math.uchicago.edu}

\thanks{The first and second authors are supported in part by NSF}

\begin{abstract}
We study the regularity of the minimizers to the problem:
\[
\la(\al, A) = \inf_{u \in H^1_0(\Om), \norm{u}_2 = 1, \abs{D} = A} \int_{\Om} \abs{Du}^2 + \al \int_{D} u^2.
\]

We prove that in the physical case $\al < \la$ in $\R^2$, any minimizer $u$ is locally $C^{1,1}$ and the boundary of the set $\set{u > c}$ is analytic where $c$ is the constant such that $D = \set{u < c}$ (up to a zero measure set).
\end{abstract}
\maketitle
\section{Introduction}\label{s-introduction}
Consider a bounded domain $\Om \subset \R^n$ with Lipschitz boundary. Fix $A$, $0 < A < \abs{\Om}$ and $\al > 0$. Our goal is to study the regularity of the minimizers to the problem:
\beq\label{e_prob_def}
\la(A, \al) = \inf_{u \in H^1_0(\Om), \norm{u}_2 = 1, \abs{D} = A} \int_{\Om} \abs{Du}^2 + \al \int_{D} u^2.
\eeq
\cite{ChanilloGIKO2000} establishes the existence of minimizers and connects \eqref{e_prob_def} with a physical problem whose goal is to minimize the first Dirichlet eigenvalue of a body of prescribed shape and mass that has to be constructed out of materials of varying densities. The Euler-Lagrange equation corresponding to \eqref{e_prob_def} is
\beq\label{e_prob_def_el}
-\lap u + \al \X_{D} u = \la(\al, A) u.
\eeq

It was proved in \cite{ChanilloGIKO2000} that for any optimal configuration $(u,D)$, there exists some $c > 0$ such that $D = \set{u < c}$ (up to a zero measure set). In fact, the weak uniqueness result in \cite{ChanilloK2007} says that this constant $c$ depends only on $\Om, \al$ and $A$, for almost every $A$.

We shall always assume here that $\al < \closure{\al}(A)$ where $\closure{\al}$ is a special constant defined in \cite{ChanilloGIKO2000}. This condition guarantees that $\al < \la(\al, A)$. The physical problem posed in \cite{ChanilloGIKO2000} in fact demands that. An elementary consequence of this condition is that $u$ is strictly superharmonic and hence satisfies the strong minimum principle. So every point in the set $\set{u = c}$ is a limit point of the set $\set{u < c}$, and $\abs{\set{u = c}} = 0$.

By a result in \cite{ChanilloGK2000}, for any point $x_0 \in \set{u = c} \cap \set{\abs{Du} > 0}$, there exists $r > 0$ such that the set $\set{u = c} \cap B_r(x_0)$ is the graph of a real-analytic function. Thus the issue is to understand points in the set $\set{u = c} \cap \set{Du = 0}$. In \cite{ChanilloK2007}, these singular points were studied for \eqref{e_prob_def_el} and a blow-up analysis performed to classify the singularities. Such an analysis was done earlier in dimension two in \cite{Blank2004} and \cite{Shahgholian2007}. The aim of this paper is to study which blow-up solutions of \cite{ChanilloK2007} are unstable for the functional \eqref{e_prob_def}. Ruling out various blow-up solutions leads therefore to improved regularity of the solution $u$ and also to regularity of the free-boundary $\set{u=c}$. In a dumb-bell shaped region $\Om$, it is proved in \cite{ChanilloGIKO2000} that one of the lobes fills faster than the other as $A \to \abs{\Om}$. Thus for certain value of $A$, one of the lobes could contain an isolated point of the set $\set{u = c}$ surrounded solely by points where $u < c$. On blow-up we will get a blow-up limit as in \cite{Shahgholian2007}, in particular the set $\set{u = c}$ is not regular. Thus in general, even if $\Om$ is simply-connected, we do not expect $\set{u = c}$ to be regular. However, it turns out that $\p \set{u > c}$ has better regularity properties. So it may be more natural to view $\p \set{u > c}$ as the free-boundary instead of $\set{u=c}$. We will therefore denote in this paper
\begin{align}
U &= \set{u > c}\\
\F &= \p U
\intertext{and}
\F^* &= \F \cap \set{\abs{Du} > 0}.
\end{align}

There is a similarity in spirit between this problem and a problem treated in \cite{MonneauW2007}. The difference being that the problem in our paper has the constraint $\abs{D} = A$, which puts complications in the construction of the variations we employ.

It will be easier to study the free functional corresponding to \eqref{e_prob_def}. We will make both variations in the domain $D$ and the function $u$. We set, for a family of domains $D(t)$ such that $\abs{D(t)} = A$,
\beq
E(s,t) = \int_{\Om} \abs{Du + s Dv}^2 + \al \int_{D(t)}(u + sv)^2 - \la \int_{\Om}(u + sv)^2.
\eeq
Our minimizing assumption then becomes
\[
E(s,t) \geq E(0,0) = 0.
\]
In section \ref{s_second_variation}, we find the formula for all first and second derivatives of $E(s,t)$. The first derivative of $E(s,t)$ with regards to $t$ already played a role in obtaining weak uniqueness in \cite{ChanilloK2007}. Pieces of the second variation formula were obtained earlier in \cite{ChanilloP}. However in order to get any contradiction the full second variation is needed.

We will confine ourselves here to state two consequences of our results. In section \ref{s_C11} we show:
\begin{thm}\label{t_intro_C11}
Let $\Om \subset \R^2$, $0 < A < \abs{\Om}$ and $0 < \al < \closure{\al}$. Let $(u,D)$ be a minimizing configuration. Then $u \in C^{1,1}(\Om)$.
\end{thm}

In contrast, one can construct solutions to the Euler-Lagrange equation \eqref{e_prob_def_el} which fail to have $C^{1,1}$ bounds \cite{ChanilloK2007}, and in a related problem, see \cite{AnderssonW2006}. We recall that \cite{Shahgholian2007} establishes that under \eqref{e_prob_def_el}, points $x_0$ where $Du(x_0) = 0$ and $U$ having positive density are isolated. Such point does exist, see \cite{ChanilloK2007}. However, we show that it is not the case for the minimizers of \eqref{e_prob_def}.

We now turn our attention to the free-boundary $\F = \p U$. We prove in section \ref{s_regularity_F} the following result:
\begin{thm}\label{t_intro_regular}
Let $(u,D)$ be a minimizing configuration. Then the set $\set{u > c}$ consists of a finite number of connected components whose closures are disjoint. The boundary of each of these components consists of finitely many disjoint, simple and closed real-analytic curves on which $\abs{Du} > 0$.
\end{thm}

The proof of theorem \ref{t_intro_regular} uses theorem \ref{t_intro_C11} and the second variation formula, but no further blow-up arguments are needed. One feature of the proof of theorem \ref{t_intro_regular} is the use of global arguments, in particular the use of the Jordan Curve Theorem. Another aspect of this problem is that one first classifies the blow-up limits and then uses the classification to get $C^{1,1}$ bounds.

It follows from these theorems and a result of \cite{ChanilloK2007} that for a minimizing configuration $(u,D)$, the 1-dim Hausdorff measure of the set $\set{u = c}$ is finite.

In the case when $\Om$ is simply connected, it follows from \cite{ChanilloGIKO2000} that $D$ is connected. From this fact and the superharmonicity of $u$, it is easy to see that each connected component of $U$ is simply connected and thus has a connected boundary. In this case, the proof of theorem \ref{t_intro_regular} simplifies considerably.

Lastly, the situation in higher dimensions is unclear. This is also the case for the problem treated in \cite{MonneauW2007}. In fact, the argument in the proof of step 2, Theorem 8.1 is incomplete because in the notation of \cite{MonneauW2007},
\[
\int_{B_1} \abs{Dw_\de}^2 \approx - \log \de \to \infty \text{ as } \de \to 0.
\]
\section{Second variation formula}\label{s_second_variation}
We start by defining what we call a regular curve. A curve $\ga: [a,b] \to \R^2$ is regular if it satisfies the following conditions
\begin{enumerate}[label=\roman*.]
\item $-\infty < a < b < \infty$
\item if $a \leq x < y \leq b$ and $x \neq a$ or $y \neq b$, then $\ga(x) \neq \ga(y)$
\item $\norm{\ga}_{C^2(a,b)}$ is finite
\item $\abs{\ga'}$ is uniformly bounded away from $0$.
\end{enumerate}

If in addition, $\ga(a) = \ga(b)$, we say it is closed and regular.

If the domain of $\ga$ is $(a,b)$, we say $\ga$ is regular (similarly closed and regular) if the continuous extension of $\ga$ to $[a,b]$ is regular (respectively closed and regular).

We state our key second variation formula in the following lemma.

\begin{lem}\label{t_second_variation}
Let $J = \cup_{k=1}^n J_k$ be a finite union of open, bounded intervals of $\R$ and
\[
\ga = (\ga_1, \ga_2): J \to \F^*
\]
a simple curve which is regular on each interval $J_k$ and $\closure{\ga(J)} \subset \F^*$. Assume also that $\dist(\ga(J_k), \ga(J_h)) > 0$ for all $1 \leq h \neq k \leq n$. For each $\xi \in J$, denote by
\[
N(\xi) = (N_1(\xi), N_2(\xi))
\]
the outward unit normal with respect to $D$ at $\ga(\xi)$. We also define $N^*$ to be $(N_2, -N_1)$ and $N'$ the first-derivative of $N$.
Let $t_0 > 0$ and $g: J \times (-t_0, t_0) \to \R$ be a function such that $g(.,t), g_t(.,t), g_{tt}(.,t) \in C(\closure{J})$ for all $t \in (-t_0, t_0)$ and
\begin{gather}
g(.,0) \equiv 0\\
\int_J g(.,t)\abs{\ga'} + \frac{1}{2}(g(.,t))^2 (N'\cdot N^*) = 0, \forall t \in (-t_0, t_0).\label{e_g_preserve_area}
\end{gather}
Then for any $v \in H^1_0$ we have
\beq\label{e_second_variation}
\left( \int_{\Om} \abs{Dv}^2 + \al \int_{D} v^2 - \La \int_{\Om} v^2 \right)
\int_{\ga} \big(g_t(\ga^{-1},0)\big)^2 \abs{Du}
\geq
\al c \left(\int_{\ga} g_t(\ga^{-1},0) v \right)^2.
\eeq
Here $g_t, g_{tt}$ denote the first and second derivatives of $g$ with respect to $t$.
\end{lem}

\begin{proof}
Reversing the direction of $\ga$ if necessary, we will assume without loss of generality that $\ga'$ and $N^*$ have the same direction, i.e
\[
\ga' \cdot N^* = \abs{\ga'}.
\]
For each $k$, it is well-known that because $\ga$ is $C^{2}$ and simple on $\closure{J_k}$, there exists a $\be_k > 0$ such that the function
\[
\phi:J_k \times [-\be_k,\be_k] \to \R^2
\]
defined by
\[
(x_1, x_2) = \phi(\xi,\be) = \ga(\xi) + \be\, N(\xi)
\]
is injective. Because $\dist(\ga(J_h), \ga(J_k)) > 0$ for $h \neq k$, we can find a number $\be_0 > 0$ such that $\phi$ is injective on $J \times [-\be_0, \be_0]$.

Substituting $t_0$ by a smaller positive number if necessary, we can assume that
\[
\norm{g}_{L^\infty(J)} < \be_0.
\]
Let
\[
K = D \setminus \phi(J\times(-\be_0,0]).
\]
Define for each $t \in (-t_0, t_0)$
\beq
D(t) = K \cup \set{\phi(\xi,\be) \st \xi \in J, \be < g(\xi,t)}
\eeq

We can compute $A(t)$, the measure of $D(t)$ by the formula
\beq\label{e_At_in_J}
A(t) = \abs{D} + \int_J \int_0^{g(\xi,t)} J(\xi,\be,t)\,d\be d\xi
\eeq
where
\begin{align*}
J(\xi,\be) &=
\begin{vmatrix}
\ga'_1 + \be  N'_1 & N_1\\
\ga'_2 + \be  N'_2 & N_2
\end{vmatrix}\\
&= \abs{ (\ga'\cdot N^*) + \be (N'\cdot N^*)}\\
&= \abs{ \abs{\ga'} + \be (N' \cdot N^*)}.
\end{align*}
Because $\norm{\ga}_{C^2(J)} < \infty$, we have $\norm{N' \cdot N^*}_{L^\infty(J)} < \infty$. Again by considering a smaller positive number $t_0$, we can assume that
\[
\norm{g}_{L^\infty(J)} \norm{N' \cdot N^*}_{L^\infty(J)} \leq \te.
\]
Thus,
\begin{align*}
\abs{\ga'} &\geq \te \geq \abs{\be}\abs{N'\cdot N^*}
\intertext{for all $\xi \in J$ and $\abs{\be} \leq \norm{g}_{L^\infty(J)}$ and so,}
J &= \abs{\ga'} + \be (N'\cdot N^*).
\end{align*}
Substituting into the formula for $A(t)$ in \eqref{e_At_in_J} we have
\begin{align*}
A(t) &= A + \int_J \int_0^{g(\xi,t)} \abs{\ga'} + \be (N' \cdot N^*)\\
&=A + \int_J g(.,t)\abs{\ga'} + \frac{1}{2} (g(.,t))^2 (N'\cdot N^*)\\
&=A \quad\text{(due to \eqref{e_g_preserve_area})}.
\end{align*}

We also have for later reference,
\begin{align*}
A'(t) &= \int_J g_t\abs{\ga'} + gg_t (N'\cdot N^*)\\
A''(t) &= \int_J g_{tt}\abs{\ga'} + (gg_{tt} + g_t^2)(N'\cdot N^*).
\end{align*}

More generally, if $F$ is a continuous function from $\R^2$ to $\R$, then
\[
\int_{D(t)}F - \int_D F = \int_J \int_0^{g(\xi,t)} F(\phi(\xi,\be)) J(\xi,\be) \,d\be d\xi
\]
and so from the Fundamental Theorem of Calculus,
\beq\label{e_dt_int_F}
\pd{}{t}\int_{D(t)}F = \int_J g_t(.,t) F(\phi(.,g(.,t))) J(.,g(.,t)).
\eeq

Define the functional
\[
E(s,t) = \int_{\Om} (Du + sDv)^2 + \al \int_{D(t)} (u+sv)^2 - \la \int_{\Om} (u + sv)^2.
\]

We will compute all second-derivatives of $E$ with respect to $s$ and $t$.

First, the second derivative of $E$ with respect to $s$,
\beq\label{e_Ess}
\pd[2]{E}{s}(s,t) = 2 \left( \int_{\Om} \abs{Dv}^2 + \al \int_{D(t)} v^2 - \la \int_{\Om} v^2 \right)
\eeq

Applying \eqref{e_dt_int_F} with $F = (u + sv)^2$, we have the first derivative of $E$ with respect to $t$,
\beq\label{e_Et}
\pd{E}{t}(s,t) = \al \int_J g_t(.,t) \big(u(\phi(.,g(.,t))) + sv(\phi(.,g(.,t)))\big)^2 J(., g(.,t))\,d\xi
\eeq

To compute the second derivative of $E$ with respect to $t$, differentiating \eqref{e_Et} and noting that
\begin{align}
\pd{}{t}\big(u(\phi(.,g(.,t)))\big)^2 &= 2u\big(\phi(.,g(.,t))\big) Du(\phi(.,g(.,t)))\cdot N g_t\notag\\
\intertext{we have}
\pd[2]{E}{t}(0,t) &= \al \pd{}{t}\int_J u(\phi(.,g))^2 g_t(\abs{\ga'} + g (N' \cdot N^*))\,d\xi\notag\\
&= \al \int_J u(\phi(.,g))^2 (g_{tt} \abs{\ga'} + (g g_{tt} + g_t^2) (N' \cdot N^*)) \notag\\
&\quad\quad {}+ \al \int_J 2u(\phi(.,g)) Du(\phi(.,g))\cdot N g_t^2 (\abs{\ga'} + g(N' \cdot N^*))\notag\\
\intertext{When $t=0$, $Du(\phi(.,g(.,0))) = Du(\ga(.)) = \abs{Du(\ga(.))} N(.)$ and so,}
\pd[2]{E}{t}(0,0) &= \al c^2 A''(0) + 2\al c \int_J g_t(.,0)^2 \abs{Du(\ga(.))}\abs{\ga'}\notag\\
&= 2\al c \int_\ga \big(g_t(\ga^{-1},0)\big)^2 \abs{Du}.\label{e_Ett}
\end{align}

To compute the mixed second derivative of $E$, differentiating \eqref{e_Et} with respect to $s$ we have
\begin{align}
\pdd{E}{s}{t}(0,t) &= 2\al \int_J g_t u(\phi(.,g)) v(\phi(.,g)) J(.,g(.,t))\notag\\
\pdd{E}{s}{t}(0,0) &= 2\al c \int_J g_t(.,0) v(\ga(.)) \abs{\ga'}\notag\\
&= 2\al c \int_\ga g_t(\ga^{-1},0) v.\label{e_Est}
\end{align}

For any value of $s$ and $t \in (-t_0, t_0)$, $u + sv \in H^1_0$ and $\abs{D(t)} = A$, so from the definition of $(u, D)$ we have that $E(0,0)$ is a minimum value of $E(s,t)$. Consequently,
\[
\pd[2]{E}{s}(0,0)\pd[2]{E}{t}(0,0) \geq \left(\pdd{E}{s}{t}(0,0)\right)^2.
\]
Substituting formula \eqref{e_Ett}, \eqref{e_Ess} and \eqref{e_Est} into this inequality we obtain the desired result.
\end{proof}

Notice that in the formula \eqref{e_second_variation}, only value of $g_t$ is present. Hence we would like to know for what kind of function $g_t$ we can find $g$ that satisfies all hypotheses of the last lemma.
\begin{lem}\label{t_second_variation_2}
Let $J$ and $\ga$ be the same as in the Lemma \ref{t_second_variation}. Assume that $h:\ga \to \R$ is a bounded, continuous function that satisfies
\[
\int_{\ga} h = 0.
\]
Then for all $v \in H^1_0(\Om)$ and $a \in \R$ we have
\begin{align}\label{e_second_variation_2}
\left( \int_{\Om} \abs{Dv}^2 + \al \int_{D} v^2 - \la \int_{\Om} v^2 \right)
\int_{\ga} h^2 \abs{Du}
&\geq
\al c \left(\int_{\ga} h (v-a) \right)^2.
\end{align}
\end{lem}

\begin{proof}
Define $N$ as in the Lemma \ref{t_second_variation}. Also define $g: J \times (-t_0, t_0) \to \R$ by
\[
g(.,t) = \frac{2 t (h\circ\ga)\abs{\ga'}}{\abs{\ga'} + \sqrt{\abs{\ga'}^2 + 2 t (h\circ\ga) \abs{\ga'} (N' \cdot N^*)}}.
\]
Since $\abs{\ga'}$ is bounded below by $\te > 0$ and $h$, $(N' \cdot N^*)$ are bounded above, we can choose $t_0$ small enough so that $g$ is well-defined in $J \times (-t_0, t_0)$. Clearly $g(.,0) \equiv 0$ and $g, g_t, g_{tt}$ are continuous functions in $J$. It also satisfies the equation
\beq\label{e_equation_of_g}
g \abs{\ga'} + \frac{1}{2}g^2 (N' \cdot N^*) = t (h \circ \ga)\abs{\ga'}
\eeq
and so for all $t \in (-t_0, t_0)$,
\[
\int_J g \abs{\ga'} + \frac{1}{2}g^2 (N' \cdot N^*) = t \int_J (h \circ \ga)\abs{\ga'} = t \int_{\ga} h = 0.
\]
Differentiating \eqref{e_equation_of_g} with respect to $t$ and letting $t=0$ we obtain
\[
g_t(.,0) = h\circ \ga.
\]
Since $g$ satisfies all the required hypothesis of the Lemma \ref{t_second_variation}, we can apply it and obtain
\begin{align*}
\left( \int_{\Om} \abs{Dv}^2 + \al \int_{D} v^2 - \la \int_{\Om} v^2 \right)
\int_{\ga} h^2 \abs{Du}
&\geq
\al c \left(\int_{\ga} h v \right)^2.
\end{align*}
Due to the fact that
\[
\int_{\ga} h = 0,
\]
we have
\[
\int_{\ga} h v = \int_{\ga} h (v-a).
\]
The conclusion then follows.
\end{proof}

\section{A regularity criterion for $\p \set{u > c}$}

\begin{lem}\label{t_regularity_criterion}
Let $P$ be a point on $\F = \p \set{u > c}$. Suppose that for each $k \in \Z^+$, there exist a positive number $r_k$, a bounded and open interval $J_k$ and a regular curve $\ga_k: J_k \to \F^*$ that satisfy the following conditions
\begin{gather*}
r_1 > r_2 > \cdots \to 0\\
\closure{\ga_k(J_k)} \subset \F^* \cap B_{r_k}(P) \setminus \closure{B_{r_{k+1}}(P)}
\end{gather*}
Then we must have
\[
\sum_{k=1}^\infty \int_{\ga(J_k)} \frac{1}{\abs{Du}} < \infty.
\]
\end{lem}

\begin{proof}
Assume without loss of generality that $P$ is the origin. Assume also that $J_k \cap J_h = \emptyset$ for all $k \neq h$ so we can use one notation $\ga$ for all $\ga_k$. We will use the following notation
\[
J_{k,m} =
\begin{cases}
J_k \cup J_{k+1} \cup \cdots \cup J_m, &\text{if }m \geq k\\
\emptyset, &\text{otherwise.}
\end{cases}
\]

Assume that
\[
\sum_{k=1}^\infty \int_{\ga(J_k)} \frac{1}{\abs{Du}} = \infty.
\]
We will derive a contradiction.

Let $V$ be a smooth, radial function in $\R^2$ such that $V$ is decreasing in $\abs{x}$ and
\beq
\begin{cases}
V(x) = 2, & \abs{x} = 0\\
2 > V(x) > 1, & \abs{x} \in (0,1/2)\\
1 > V(x) > 0, & \abs{x} \in (1/2,1)\\
V(x) = 0, & \abs{x} \geq 1\\
\end{cases}
\eeq
For each $k \in \Z^+$, define $v_k(x) = V(x/r_k)$. It is easy to verify that when $r_k$ is small enough,
\[
\int_\Om \abs{Dv_k}^2 = \int_\Om \abs{DV}^2
\]
and so for any $k$ large enough
\begin{align}
\int_{\Om} \abs{Dv_k}^2 + \al \int_D \abs{v_k}^2 - \la \int_{\Om} \abs{v_k}^2 &< \int_\Om \abs{DV}^2 < \infty\label{e_uniform_bound_Ev}.
\end{align}
We will drop the subscript $k$ from the rest of the proof. We list here values of $v-1$ for easy reference later,
\beq
\begin{cases}
v(x) - 1 = 1, & \abs{x} = 0\\
1 > v(x) - 1 > 0, & \abs{x} \in (0, r_k/2)\\
0 > v(x) - 1 > -1, & \abs{x} \in (r_k/2, r_k)\\
v(x) - 1 = -1, & \abs{x} \geq r_k.\\
\end{cases}
\eeq

Because $J_k$ and $\abs{\ga'}$ are bounded, $\ga(J_k)$ is of finite length. We also have $\abs{Du}$ is uniformly bounded away from 0 on $\ga(J)$ since $\closure{\ga(J)} \subset \F^*$. Together with the fact that $\ga(J_{0,k-1}) \subset \compl{B_{r_k}}$, we have
\[
-\infty < \int_{\ga(J_{0,k-1})} \frac{v - 1}{\abs{Du}} = -\int_{\ga(J_{0,k-1})} \frac{1}{\abs{Du}} < 0.
\]
Choose an $m$ such that $r_m < r_k /2$. From the facts that $v(x) - 1 > 0 \text{ in } B_{r_m}$, $\ga(J_l) \subset B_{r_m}$ for all $l \geq m$ and $v(x) - 1 \to 1$ as $\abs{x} \to 0$ we have
\[
\int_{\ga(J_{m,\infty})} \frac{v - 1}{\abs{Du}} \sim \int_{\ga(J_{m, \infty})} \frac{1}{\abs{Du}} = \infty.
\]

Consequently, there must be a number $l \geq m$ such that
\[
\int_{\ga(J_{m,l-1})} \frac{v-1}{\abs{Du}} \leq
-\int_{\ga(J_{0,k-1})} \frac{v - 1}{\abs{Du}} <
\int_{\ga(J_{m,l})} \frac{v-1}{\abs{Du}}.
\]
Choose a subinterval $J'_l \subset J_l$ such that
\[
\int_{\ga(J_{m,l-1})} \frac{v-1}{\abs{Du}}
+ \int_{\ga(J'_l)} \frac{v - 1}{\abs{Du}} =
- \int_{\ga(J_{0,k-1})} \frac{v - 1}{\abs{Du}}.
\]
In other words, we have
\[
\int_{\ga(J^k)} \frac{v - 1}{\abs{Du}} = 0.
\]
where $J^k = J_{0,k-1} \cup J_{m,l-1} \cup J'_l$.

We can now apply the Lemma \ref{t_second_variation_2} to $J^k$, $\ga$, $v$, $a=1$ and $h = (v-1)/\abs{Du}$, and obtain
\begin{align*}
\int_\Om \abs{DV}^2 \int_{\ga(J^k)} \frac{(v - 1)^2}{\abs{Du}} &\geq
 \al c\left(\int_{\ga(J^k)} \frac{(v - 1)^2}{\abs{Du}}\right)^2\\
\int_\Om \abs{DV}^2 &\geq \al c\int_{\ga(J^k)} \frac{(v - 1)^2}{\abs{Du}}\\
&\geq \al c\int_{\ga(J_{0,k-1})} \frac{(v - 1)^2}{\abs{Du}}\\
&\geq \al c\int_{\ga(J_{0,k-1})} \frac{1}{\abs{Du}} \quad\text{($v-1=-1$ on $\ga(J_{0,k-1}) \subset \compl{B_{r_k}}$)}.\\
\intertext{Let $k$ go to $\infty$ we have}
\int_\Om \abs{DV}^2 &\geq \al c\int_{\ga(J_{0,\infty})} \frac{1}{\abs{Du}}
= \infty
\end{align*}
which is a contradiction.

So we must have
\[
\sum_{k=1}^\infty \int_{\ga(J_k)} \frac{1}{\abs{Du}} < \infty
\]
as desired.
\end{proof}


Next, we prove a direct consequence of the last lemma. Informally, it says that if the set $\p \set{u > c} \cap \set{\abs{Du} > 0}$ is big enough around a point of $\p \set{u > c}$, then at this point, $\abs{Du} > 0$.
\begin{lem}\label{t_regularity_criterion_2}
Let $P$ be a point on $\F = \p \set{u > c}$. Suppose that there are numbers $K \in \Z$ and $\si > 0$ such that for each $k \geq K$, there exists a regular curve $\ga_k:J_k \to \F^*$ with the following properties
\begin{gather*}
\closure{\ga_k(J_k)} \subset \F^* \cap B_{2^{-k}}(P) \setminus \closure{B_{2^{-(k+1)}}(P)}\\
\hm{1}(\ga_k(J_k)) = \int_{J_k} \abs{\ga'_k} > \si 2^{-k}.
\end{gather*}
Then $\abs{Du(P)} > 0$.
\end{lem}

\begin{proof}
Assume that $Du(P) = 0$. To derive a contradiction, it is enough to show that
\[
\sum_{k=K}^\infty \int_{\ga(J_k)} \frac{1}{\abs{Du}} = \infty
\]
and use the last lemma.

From a result in \cite{Chemin1998} and the fact that $\lap u \in L^\infty$, there exists some positive constant $C$ such that for all $x \in \Om$,
\[
\abs{Du(x)} = \abs{Du(x) - Du(P)} \leq C \abs{x - P} \log (1 / \abs{x - P}).
\]
Thus,
\begin{align*}
\sum_{k=K}^\infty \int_{\ga(J_k)} \frac{1}{\abs{Du}} &\geq \frac{1}{C} \sum_{k=K}^\infty \int_{\ga(J_k)} \frac{1}{\abs{x - P} \log(1/\abs{x - P})}\\
&\geq \frac{1}{C} \sum_{k=K}^\infty \frac{\si 2^{-k}}{2^{-k}\log (2^k)}\\
&= \frac{1}{C} \frac{\si}{\log 2} \sum_{k=K}^\infty \frac{1}{k}\\
&= \infty.
\end{align*}
\end{proof}

\section{$C^{1,1}$ regularity}\label{s_C11}

We now apply the regularity criterion from the last section to show that if the set $\set{u > c}$ has positive density at a point of the set $\p \set{u > c}$, then at that point $\abs{Du} > 0$.
\begin{thm}\label{t_regularity_positive_density}
Let $P$ be a point on $\p \set{u > c}$. Assume that there exist $\be, r_0 > 0$ such that
\[
\abs{\set{u>c} \cap B_r(P)} \geq \be r^2
\]
for all $0 < r < r_0$. Then $\abs{Du(P)} > 0$.
\end{thm}

\begin{proof}
Without loss of generality, let $P$ be the origin. Assume that $Du(P) = 0$. For each $r > 0$, define
\[
v_r(x) = \frac{c - u(rx)}{r^2}.
\]
Also define $I(r)$ to be the supremum of lengths of all regular curves with closures in the set
\[
\set{v_r = 0} \cap \set{\abs{Dv_r} > 0} \cap B_1 \setminus \closure{B_{1/2}}.
\]
We show that there exist some $r_0 > 0$ and $\si > 0$ such that $I(r) > \si$ for all $0 < r < r_0$.

Assume that it is not the case, then there exists a sequence $r_k \to 0$ such that $I(r_k) \to 0$. As a consequence of the Theorem 3.1 in \cite{ChanilloK2007}, two possibilities arise.
\begin{enumerate}
\item A subsequence of $v_{r_k}/T(r_k)$ converges to a non-zero, homogeneous of degree 2 harmonic function where
\[
T(r) = \frac{1}{r^2}\left(\frac{1}{2\pi r}\int_{\p B_r} (c - u)^2\right)^{1/2}.
\]
\item A subsequence of $v_{r_k}$ converges to a homogeneous solution of degree 2 of the equation
\[
\lap v = c(\la - \al) \X_{\set{v \geq 0}} + c\la \X_{\set{v < 0}}.
\]
\end{enumerate}

We consider case (1) first. Without loss of generality, we can assume that $v_{r_k}/T(r_k)$ converges to $v(x) = x_1 x_2$ in $C^{1,\de}$ as $k \to \infty$. We will hereafter denote $v_{r_k}$ by $v_k$ and $T(r_k)$ by $T_k$.

Let $\ep$ be any number in $(0,1/8)$. It can be verified easily that
\[
Q_1 = [1/2 + \ep, 1 - \ep] \times [-\ep, \ep]
\]
is a subset of the set $B_1 \setminus \closure{B_{1/2}}$. We have for any $x_1 \in [1/2 + \ep, 1 - \ep]$,
\begin{enumerate}[label=\roman*.]
\item The first-derivative with respect to $x_2$, $v_2(x_1,.) = x_1 \in (1/2, 1)$.
\item $v(x_1,-\ep) = -\ep x_1 \leq -\ep/2$ and $v(x_1,\ep) = \ep x_1 \geq \ep/2$.
\end{enumerate}
Since $v_k/T_k \to v$ in $C^{1,\de}$, we can choose some $N$ such that for all $k > N$,
\[
\abs{v_k/T_k - v}_{\infty} < \ep/4, \abs{(v_k)_2/T_k - v_2} < 1/4 \text{ in $Q_1$}.
\]
It follows that for all $k > N$,
\begin{enumerate}
\item $5/4 > (v_k)_2/T_k > 1/4$ on $[1/2+\ep, 1 - \ep] \times [-\ep, \ep]$.
\item $v_k(x_1,-\ep)/T_k \leq -\ep/4$ and $v_k(x_1,\ep)/T_k \geq \ep/4$.
\end{enumerate}
Consequently, for each $x_1$, there is exactly one value of $x_2$ such that $v_k(x_1,x_2) = 0$. Denote this value by $\tau_k(x_1)$ and define $\ga_k(x_1) = (x_1, \tau_k(x_1))$. Since $5/4 > (v_k)_2 / T_k > 1/4$, doing implicit differentiation we have $-\infty < \tau_k' < \infty$ and so $1 \leq \abs{\ga_k'} < \infty$. $\ga_k$ is also clearly the boundary of a connected component of the set $\set{v_k < 0}$ since a neighborhood below it is an open subset of the set $\set{v_k < 0}$. The length of $\ga_k$ is at least
\[
(1-\ep) - (1/2 + \ep) = 1/2 - 2\ep > 1/4.
\]
This implies that $I(r_k) > 1/4$ for all $k > N$, contradicting our assumption that $I(r_k) \to 0$.

In the second case, we can also assume that $v_k$ converges in $C^{1,\de}$ to a homogeneous solution of degree 2 of the equation
\[
\lap v = c(\la - \al) \X_{\set{v \geq 0}} + c\la \X_{\set{v < 0}}.
\]
Since
\[
\abs{\set{u>c} \cap B_r} \geq \be r^2,
\]
in terms of $v_k$ we have
\[
\abs{\set{v_k < 0} \cap B_1} \geq \be.
\]
Letting $k$ go to $\infty$ we obtain
\[
\abs{\set{v \leq 0} \cap B_1} \geq \be.
\]
From the Lemma 1.2 in \cite{Shahgholian2007}, we know that either the set $\set{v=0} \cap \set{Dv = 0} = \set{0}$ or $v$ is of the form $c(\la - \al) x_1^2/2$ after a rotation. Because
\[
\abs{\set{c(\la - \al) x_1^2/2 \leq 0} \cap B_1} = 0,
\]
contradicting the positive density condition for $v$ above, we must have then
\[
\set{v=0} \cap \set{Dv = 0} = \set{0}.
\]

Since
\[
\abs{\set{v \leq 0} \cap B_1} \geq \be,
\]
$v$ is superharmonic and $v$ is homogeneous, there exists a point $z$ such that $\abs{z} = 1$ and $v(tz) = 0$ for all $t \in [0,1]$. Assume that $z = (1,0)$. We also have $Dv(1/2,0) \neq 0$ due to the fact that
\[
\set{v=0} \cap \set{Dv = 0} = \set{0}.
\]
Because $\abs{v_2(1/2,0)} = \abs{Dv(1/2,0)} \neq 0$, we can assume without loss of generality that $u_2(1/2,0) > 0$. Now, arguing similarly to the first case, we obtain $I(r_n) > \si$ for some $\si > 0$ when $n$ large enough, contradicting our assumption that $I(r_n) \to 0$.

Thus, in all cases, there exist $\si > 0$ and $r_0 > 0$ such that $I(r) > \si$ for all $0 < r < r_0$. In other words, for each $r < r_0$, there exists a regular curve of length at least $\si$ with closure in the set
\[
\set{v_r = 0} \cap \set{Dv_r \neq 0} \cap B_1 \setminus \closure{B_{1/2}}.
\]
In terms of $u$, it means for all $0 < r < r_0$, there exists a regular curve of length at least $\si r$ with closure in the set
\[
\F^* \cap B_r \setminus \closure{B_{r/2}}.
\]
Applying the Lemma \ref{t_regularity_criterion_2} we have $\abs{Du(P)} > 0$, contradicting the assumption that $Du(P) = 0$. So $\abs{Du(P)} > 0$.
\end{proof}


\begin{cor}\label{t_C11}
$u \in C^{1,1}(\Om)$.
\end{cor}

\begin{proof}
It is clear that $u$ is $C^{1,1}$ at points in $\set{u \neq c}$ or $\set{u = c} \cap \set{\abs{Du} > 0}$. Assume that there exists a point $P \in \set{u = c} \cap \set{\abs{Du} = 0}$ at which $u$ is not $C^{1,1}$. In other words,
\[
\limsup_{r \to 0} \sup_{\abs{x-P} < r} \frac{\abs{u(rx + P) - c}}{r^2} \to \infty.
\]
From the Lemma 3.18 in \cite{ChanilloK2007}, there must exist $\be, r_0 > 0$ such that
\[
\abs{\set{u > c} \cap B_r(P)} \geq \be r^2 \text{ for all } 0 < r < r_0.
\]
However, the Theorem \ref{t_regularity_positive_density} then implies that $\abs{Du(P)} > 0$, a contradiction. Thus $u \in C^{1,1}(\Om)$.
\end{proof}

\section{Regularity of connected components of $\set{u > c}$}\label{s_regularity_component_U}

We first prove the following lemma.
\begin{lem}\label{t_locally_regular_curve}
Let $L \subset \R^2$ be a connected set. Furthermore, assume that for any $P \in L$, there exists $r > 0$ such that the set $B_r(P) \cap L$ is a regular curve. Then given any pair $S, Q \in L$, there exists a regular curve in $L$ with $S, Q$ as two end points.
\end{lem}

\begin{proof}
Define $L_S$ to be the set of points $R \in L$ such that there exists a regular curve in $L$ with $S, R$ as two endpoints. We will show that $L_S$ is non-empty, closed and open. Because $L$ is connected, it means $L_S = L$ and the conclusion follows.

Let $r > 0$ be a number such that $L \cap B_r(S)$ is a regular curve. Obviously, any point in this set is a point in the set $L_S$ as well. So $L_S$ is non-empty.

Assume that $R \in L_S$. Let $r > 0$ be a number such that $B_r(R) \cap L$ is a regular curve. Because there is a regular curve connecting $S$ and $R$, it is easy to see that for any $R' \in B_r(R) \cap L$, we can truncate or extend that regular curve to obtain a new regular curve connecting $S$ and $R'$. Thus, $B_r(R) \cap L \subset L_R$. Since it is true for all $R \in L_S$, $L_S$ must be open.

Arguing similarly we have, if $R \in \compl{L_S}$, then there exists $r > 0$ such that $B_r(R) \cap L \subset \compl{L_S}$. In other words, $L_S$ is closed.
\end{proof}

Next, we prove our first result about the structure of the set $\p \set{u > c} \cap \set{\abs{Du} > 0}$.
\begin{thm}\label{t_regularity_whole_component}
If $\F_1$ is a connected component of $\F = \p U$, then either $\abs{Du} > 0$ at every point of $\F_1$, or $\abs{Du} \equiv 0$ on $\F_1$.
\end{thm}

\begin{proof}
Assume that $\F_1$ contains at least one point where $\abs{Du} > 0$. Let $L$ be a connected component of the set $\F_1 \cap \set{\abs{Du} > 0}$. $L$ must be non-empty by definition.

Since for each $S \in L$, there exists a number $r > 0$ such that $B_r(S) \cap \p \set{u > c}$ is a simple, analytic curve where $\abs{Du} > 0$, $L$ has to be open.

We will show that $L$ is closed as well. Choose any convergent sequence $\set{P_n}$ in $L$. Because $\F_1$ is a connected component of $\F$, $\F_1$ is closed. Thus, there exists some $P \in \F_1$ such that
\[
P_n \to P \in \F_1 \text{ as } n \to \infty.
\]
Pick any $r_0 < \abs{P_1 - P}$ (here $P_1$ is the first point in the sequence $\set{P_n}$). For any $0< r < r_0$, there exists some $P_n$ such that $\abs{P_n - P} < r/2$. From Lemma \ref{t_locally_regular_curve} we have there exists a regular curve $\ga:[0,l] \to L$ such that $\ga(0) = P_1$ and $\ga(l) = P_n$. Define
\begin{align*}
a &= \inf \set{s \in [0,l] \st \ga([s,l]) \subset B_r(P)}\\
b &= \inf \set{s \in [a,l] \st \abs{\ga(s) - P} = r/2}.
\end{align*}
The existence of $a < b \in (0,l)$ is justified because $\ga$ is a regular curve and $\abs{\ga(0) - P} > r$ while $\abs{\ga(l) - P} < r/2$. It can also be verified easily that
\begin{gather*}
\abs{\ga(a) - P} = r, \abs{\ga(b) - P} = r/2\\
\ga((a,b)) \subset B_r(P) \setminus \closure{B_{r/2}(P)}\\
\hm{1}(\ga((a,b))) \geq r/2.
\end{gather*}
Pick some $\ep > 0$ small so that the length of the segment $\ga((a+\ep, b - \ep))$ is at least $r/3$. It also follows from the above argument that
\[
\closure{\ga((a+\ep, b - \ep))} \subset \F^* \cap B_r(P) \setminus \closure{B_{r/2}(P)}.
\]
Since we can do it for all $r < r_0$, the Lemma \ref{t_regularity_criterion_2} then implies that $\abs{Du(P)} > 0$. Consequently, there exists $r_1 > 0$ such that $B_{r_1}(P) \cap \F$ is a regular curve where $\abs{Du} > 0$. Since $\F_1$ is a connected component of $\F$ and $P \in \F_1$, the whole curve $B_{r_1}(P) \cap \F$ must be in $\F_1$. Pick some $P_n$ such that $\abs{P_n - P} < r_1$. It is clear that $P_n$ has to be in the curve $B_{r_1}(P) \cap \F$. But because $P_n \in L$, $\abs{Du} > 0$ on $B_{r_1}(P) \cap \F$ and $L$ is connected, the whole curve $B_{r_1}(P) \cap \F$ has to be in $L$. In particular, $P \in L$. Since $\set{P_n}$ is an arbitrary convergent sequence in $L$, it implies that $L$ is closed.

We have proved that $L$ is non-empty, open and closed. Because $\F_1$ is connected, we have $L = \F_1$. In other words $\abs{Du} > 0$ for every point on $\F_1$.
\end{proof}


\begin{lem}\label{t_component_pU1_F}
Let $U_1$ be a connected component of $U$ and $\F_1$ a connected component of $\p U_1$ such that $\abs{Du} > 0$ on $\F_1$. Then $\F_1$ is also a connected component of $\F$.
\end{lem}

\begin{proof}
Let $P$ be any point on $\F_1$. Because $\abs{Du(P)} > 0$, there exists $r > 0$ such that the set $B_r(P) \cap \F$ is a regular curve that divides $B_r(P)$ into two disjoint connected regions, one where $u < c$ and one where $u > c$. It is easy to see that the connected region where $u > c$ is a subset of $U_1$ and so $B_r(P) \cap \F \subset \F_1$. Now for each point in $\F_1$, pick a ball like before and consider the union $V$ of all these balls. Clearly $V$ is open and $V \cap \F = \F_1$. Hence, $\F_1$ is a connected component of $\F$.
\end{proof}


Next we show that the set $\set{\abs{Du} > 0}$ is dense in the boundary of each connected component of $U$, improving Lemma 2.3 in \cite{ChanilloK2007}.

\begin{lem}\label{t_F*_dense_pU1}
If $U_1$ is a connected component of $U$, then
$\p U_1 = \closure{\p U_1 \cap \set{\abs{Du} > 0}}$.
\end{lem}

\begin{proof}
Let $P$ be a point on $\p U_1$ such that $Du(P) = 0$. We will show that for any $\ep > 0$, there exists a point $Q \in \p U_1$ such that $\abs{P - Q} < \ep$ and $\abs{Du(Q)} > 0$.

Since $P \in \p U_1$, we can choose a point $S \in U_1$ such that $\abs{P - S} < \ep/2$. Define
\[
r = \sup \set{s \st B_s(S) \subset U_1}.
\]
It is obvious that $0 < r \leq \abs{P - S} < \ep/2$ and $\p B_r(S) \cap \p U_1 \neq \emptyset$. Let $Q$ be any point of the set $\p B_r(S) \cap \p U_1$. Because $u$ is superharmonic and $Q$ is a boundary minimum point of $u$ in the set $\closure{B_r(z)}$, from the Hopf's Lemma we have $\abs{Du(Q)} > 0$. We also have easily $\abs{P - Q} < \ep$ due to the facts that $\abs{P - S} < \ep/2$ and $\abs{S - Q} = r < \ep/2$.
\end{proof}


\begin{lem}\label{t_small_component_pU1}
Let $U_1$ be a connected component of $U$ and $P$ a point on $\p U_1$ such that $\abs{Du(P)} = 0$, then for any $r > 0$, there exists a connected component $\F_1$ of $\p U_1$ such that $\abs{Du} > 0$ in $\F_1$ and $\F_1 \subset B_r(P)$.
\end{lem}

\begin{proof}
Let's assume that $P$ is the origin. First, we show that there exists an $r' > 0$ such that for any connected component $\F_1$ of $\F$ where $\abs{Du} > 0$, if $\F_1 \cap \compl{B_r} \neq \emptyset$, then $\F_1 \subset \compl{B_{r'}}$. In other words, if $\F_1$ contains a point outside $B_r$, then the whole component $\F_1$ has to stay outside $B_{r'}$.

If it is not the case, then for any $r' > 0$, there exists some connected component $\F_1$ of $\F$ such that $\abs{Du} > 0$ on $\F_1$, $\F_1 \cap \compl{B_r} \neq \emptyset$ and $\F_1 \cap B_{r'} \neq \emptyset$. It means for any $k > \log_2 (1/r)$, there exists a connected component $\F_1$ of $\F$ such that $\F_1 \cap \compl{B_{2^{-k}}} \neq \emptyset$, $\F_1 \cap B_{2^{-k-1}} \neq \emptyset$ and $\abs{Du} > 0$ on $\F_1$. Choose $P_1, P_2 \in \F_1$ such that $\abs{P_1} \geq 2^{-k}$, $\abs{P_2} < 2^{-(k+1)}$. From the Lemma \ref{t_locally_regular_curve}, there exists a regular curve connecting $P_1$ and $P_2$. Arguing as in the Lemma \ref{t_regularity_whole_component}, we can find a smaller regular piece of this curve of length at least $2^{-k}/3$ in the set
\[
\F_1 \cap B_{2^{-k}} \setminus \closure{B_{2^{-k-1}}}.
\]
Since we can do it for all $k > \log_2 (1/r)$, applying the Lemma \ref{t_regularity_criterion_2} we can conclude that $\abs{Du(P)} > 0$, contradicting our hypothesis on $P$. The existence of $r'$ then follows.

Now using the Lemma \ref{t_F*_dense_pU1}, we can choose a point $Q \in \p U_1$ such that $Q \in B_{r'}$ and $\abs{Du(Q)} > 0$. Let $\F_1$ be the connected component of $\p U_1$ that contains $Q$. It follows from what we just proved above that $\F_1 \subset B_r$. To show that $\abs{Du} > 0$ on $\F_1$, just note that because $\F_1$ is a connected component of $\p U_1$ and $\p U_1 \subset \F$, there exists a connected component $\F'_1$ of $\F$ such that $\F_1 \subset \F'_1$. Because $\abs{Du(Q)} > 0$ and $Q \in \F_1 \subset \F'_1$, applying the Lemma \ref{t_regularity_whole_component} we have $\abs{Du} > 0$ on $\F'_1$.
\end{proof}


Next, we prove a lemma about the geometric structure of regular connected components of $\F$.
\begin{lem}\label{t_regular_component_closed_curve}
If $\F_1$ is a connected component of $\F$ such that $\abs{Du} > 0$ on $\F_1$, then $\F_1$ is a closed and regular curve.
\end{lem}

\begin{proof}
Pick any point $P$ on $\F_1$. Consider the ODE
\beq
\ga'(t) = \frac{(Du(\ga(t)))^*}{\abs{Du(\ga(t))}}, \ga(0) = P
\eeq
where $\ga$ is a function from $[0,\infty)$ to $\F_1$. Here, as in the section \ref{s_second_variation}, $N^*$ denotes the vector obtained from rotating $N$ clockwise an angle of $\pi/2$.

First, it is easy to see that if a solution $\ga$ exists up to some time $t_0$, then we can extend that solution to $t_0 + \ep$ for some $\ep > 0$. Indeed, because $\F_1$ is closed, so $\ga(t_0) \subset \F_1$. Since $\F_1$ is regular, there exists some $r > 0$ such that $\F_1 \cap B_r(\ga(t_0))$ is the graph of a analytic function. Thus, we can extend $\ga$ to some time $t_0 + \ep$. Consequently, this solution $\ga$ exists for all time.

Define
\[
T = \sup \set{t \st \ga((0,t)) \text{ is simple}}.
\]
Because $B_r(P) \cap \F_1$ is a simple curve for some $r > 0$ small, $T \geq r > 0$. We also have since $\abs{\ga'} = 1$ that the length of $\ga((0,T))$ is exactly $T$. We will show $T < \infty$ by proving that $\hm{1}(\F_1) < \infty$.

Since $\F_1 \subset \F^*$, for each point $Q \in \F_1$, there exists $r > 0$ such that $B_r(Q) \cap \F_1$ is an analytic curve. It implies that $\hm{1}(B_r(Q) \cap \F_1) < \infty$. Because $\F_1$ is closed and bounded, we can cover $\F_1$ by a finite number of such balls and so $\hm{1}(\F_1) < \infty$.

We will show that there exists a time $T' \in [0,T)$ such that $\ga(T') = \ga(T)$.

Choose a decreasing sequence of $\set{t_k}$ that converges to $T$. Define
\begin{align*}
a_k &= \inf \set{a \in [0,t_k) \st \ga(a) = \ga(t) \text{ for some } t \in (a,t_k)}
\intertext{and}
b_k &= \inf \set{b \in (a_k,t_k) \st \ga(a_k) = \ga(b)}.
\end{align*}

The existence of $a_k$ is justified from the fact that $\ga([0,t_k))$ is not simple. The existence of $b_k \geq a_k$ follows the continuity of $\ga$. We show that actually $b_k > a_k$. Indeed, since there exists an $r > 0$ such that $B_r(\ga(a_k)) \cap \F_1$ is a simple curve, there is no $t \in (a_k, a_k + r)$ such that $\ga(a_k) = \ga(t)$ and so $b_k \geq a_k + r > a_k$.

We also have other properties of $a_k, b_k$
\begin{enumerate}[label=\roman*.]
\item $\set{a_k}$ is increasing.
\item $a_k \leq T \leq b_k < t_k$.
\end{enumerate}
Passing to a subsequence if necessary, assume that $a_k \to T'$ as $k \to \infty$. It is trivial that $b_k \to T$ and $\ga(T') = \ga(T)$. All we need to do now is to show that $T' < T$. Indeed since there exists $r > 0$ such that $B_r(\ga(T)) \cap \F_1$ is a simple curve, $\ga((T-r,T+r))$ is a simple curve. When $k$ is large enough, $b_k \in [T,T+r)$ and consequently $a_k \leq T-r$. Thus $T' \leq T - r < T$. We also note that there exists no other pair $(a,b) \neq (T', T)$ with $0 \leq a < b \leq T$ such that $\ga(a) = \ga(b)$.

If $T' \neq 0$, then as a consequence of the result above, for all $r > 0$ small, the set $\F_1 \cap B_r(\ga(T'))$ consists of three disjoint arcs $\ga(T' - r, T']$, $\ga[T', T' + r)$ and $\ga(T-r, T]$ that intersect at an endpoint $\ga(T')$, contradicting the fact that $B_r(\ga(T')) \cap \F_1$ is a regular curve when $r > 0$ is small. Thus, $T' = 0$.

To show that $\F_1 = \ga([0,T])$, we argue the same way as in the Lemma \ref{t_component_pU1_F} to show that there exists an open set $V$ such that $V \cap \ga([0,T]) = \ga([0,T])$ and note that $\F_1$ is connected.
\end{proof}

If $\F_1$ is a connected component of $\F$ such that $\abs{Du} > 0$, then by the Lemma \ref{t_regular_component_closed_curve} above, we know that $\F_1$ is a closed and regular curve. Using the Jordan Curve Theorem (see for example \cite{Whyburn1964}), we know that $\F_1$ divides $\R^2$ into two separate regions, an inside region and an outside region. We will denote the inside region as $I(\F_1)$ and the outside region $O(\F_1)$.


\begin{thm}
If $U_1$ is a connected component of $U$, then $\abs{Du} > 0$ on $\p U_1$.
\end{thm}

\begin{proof}
Without loss of generality, let's assume that $0 \in \p U_1$ and $Du(0) = 0$. Choose $r > 0$ such that
\[
r < \frac{c}{\norm{Du}_\infty} \text{ and } U_1 \not\subset B_r
\]

From the Lemma \ref{t_small_component_pU1} we have that there exists some connected component $\F_1$ of $\p U_1$ such that $\abs{Du} > 0$ in $\F_1$ and $\F_1 \subset B_r$. By the Lemmas \ref{t_component_pU1_F} and \ref{t_regular_component_closed_curve}, $\F_1$ is closed and regular. Thus, following the remark preceding this theorem, we can talk about the inside region $I(\F_1)$ and outside region $O(\F_1)$. We know that both $I(\F_1)$ and $O(\F_1)$ are open and connected. Furthermore, $I(\F_1)$ is bounded while $O(\F_1)$ is unbounded. Because $\F_1 \subset B_r$, we can connect any point in $\compl{B_r}$ to a point far away by a line that does not intersect $\F_1$ and so $\compl{B_r} \subset O(\F_1)$. Consequently, $I(\F_1) \subset B_r$.

Because $U_1$ is connected, we must have either $U_1 \subset I(\F_1)$ or $U_1 \subset O(\F_1)$. Since $I(\F_1) \subset B_r$ and $U_1 \not\subset B_r$, we cannot have $U_1 \subset I(\F_1)$. Thus $U_1 \subset O(\F_1)$. Let $P$ be a point on $\F_1$. There exists $r' > 0$ such that $B_{r'}(P) \cap \F_1$ is a regular curve that divides $B_{r'}(P)$ into two disjoint connected regions, one where $u > c$ and another where $u < c$. Because $P$ is a boundary point of $U_1$, it is clear that the region where $u > c$ must be a subset of $U_1$ and so, a subset of $O(\F_1)$. It implies that the region where $u < c$ is a subset of $I(\F_1)$. Thus $u < c$ for some point in $I(\F_1)$. However, since $u$ is superharmonic, $u$ cannot have an interior minimum in the set $I(\F_1)$. Thus, there must be a point $Q \in I(\F_1)$ such that $u(Q) = 0$. In other words, $Q \in \p \Om$. But then from the facts that $Q \in I(\F_1) \subset B_r$ and
\[
r < \frac{c}{\norm{Du}_\infty}
\]
we must have
\[
\abs{u(Q) - u(0)} < r\abs{Du}_\infty < c,
\]
contradicting the fact that $u(Q) = 0$ and $u(0) = c$.

In other words, $\abs{Du} > 0$ at every point on $\p U_1$.
\end{proof}

\section{Regularity of $\p \set{u > c}$}\label{s_regularity_F}
At the end of last section, we have proved that $\abs{Du} > 0$ on the boundary of each component of $U$. It might still happen that connected components of $U$ accumulate to a point where $\abs{Du} = 0$. For example, connected components of $U$ consists a sequence of smaller and smaller balls that converge to a point. In this section, we prove that this scenario cannot happen. Indeed, $U$ only has a finite number of connected components.


\begin{lem}
Let $U_1$ be a connected component of $U$. Then there exists a unique connected component $\F_1$ of $\p U_1$ such that $U_1 \subset I(\F_1)$. We will say that $\F_1$ surrounds $U_1$.
\end{lem}

\begin{proof}
Pick any point $P \in U_1$. Define
\[
d = \sup \set{\abs{P - x} \st x \in U_1}.
\]
Clearly, there exists a point $Q \in \p U_1$ such that $\abs{P-Q} = d$. Assume without loss of generality that $P$ is the origin and $Q = (d, 0)$. It is easy to see that $U_1$ has to be on the left-side of the line $x_1 = d$ due to the definition of $d$. From this and the fact that $(d,0) \in \p U_1$, we have the outward unit normal with respect to $U_1$ at $Q$ has to be $e_1$. Let $\F_1$ be the connected component of $\p U_1$ that contains $Q$. Note that $e_1$ will also be the outward unit normal to $I(\F_1)$ and so, there must exist some $\ep > 0$ such that $(d, d - \ep)\times \set{0} \subset I(\F_1)$ and $(d + \ep, d)\times\set{0} \subset O(\F_1)$.

Let $r > 0$ such that $B_r(Q) \cap \F_1$ is a regular curve that divides $B_r(Q)$ into two regions, $u > c$ and $u < c$. Since $U_1$ is connected and $Q \in \p U_1$, the region $u > c$ is a subset of $U_1$. Because the outward unit normal vector at $Q$ to this curve is $e_1$, by choosing a smaller $\ep$ if necessary, we have $u < c$ on one of two sets $(d, d - \ep)\times \set{0}$, $(d, d + \ep)\times \set{0}$ and $u > c$ on the other. Because the set where $u > c$ must be a subset of $U_1$, it has to be on the left-side of $(d,0)$ and thus, it has to be $(d, d - \ep)\times \set{0}$. Hence $U_1 \cap I(\F_1) \neq \emptyset$. But $U_1$ is connected, so $U_1 \subset I(\F_1)$.

Assume there is another connected component $\F_2$ of $\p U_1$ such that $U_1 \subset I(\F_2)$. It is easy to derive that $\F_1 \subset \closure{I(\F_2)}$ and $\F_2 \subset \closure{I(\F_1)}$. Consequently, $\F_2 \equiv \F_1$.
\end{proof}


\begin{lem}\label{t_lower_bound_int_1_Du}
Let $U_1$ be a connected component of $U$ and $\F_1$ the connected component of $\p U_1$ that surrounds $U_1$. Assume further that $u \geq c/2$ in the convex hull of $I(\F_1)$. Then
\[
\int_{\F_1} \frac{1}{\abs{Du}} \geq \frac{1}{C_1}
\]
where $C_1 = \norm{u}_{C^{1,1}(\set{u \geq c/2})}$.
\end{lem}

\begin{proof}
Without loss of generality, assume that $u$ attains its maximum value in $U_1$ at the origin. Let $P$ be the point on $\F_1$ such that
\[
\abs{P} = \max \set{\abs{x} \st x \in \F_1}.
\]
Let $x$ be any point on $\F_1$. Since both $x$ and $0$ belongs to the convex hull of $I(\F_1)$, $u \geq c/2$ on the line segment that connects $0$ and $x$. Thus, we have
\begin{align*}
\abs{Du(x)} &= \abs{Du(x) - Du(0)}\\
&\leq C_1\abs{x}\\
&\leq C_1\abs{P}.
\end{align*}
Because $P \in \F_1$ and $0 \in I(\F_1)$, the line connecting $P$ and $0$ has to intersect with $\F_1$ at another point $Q$ and $0$ is between $P$ and $Q$. Clearly, the length of $\F_1$ is greater than the length of the line segment $PQ$ which is greater than $\abs{P}$. Thus,
\[
\int_{\F_1}\frac{1}{\abs{Du}} > \frac{\abs{P}}{C_1\abs{P}} = \frac{1}{C_1}.
\]
\end{proof}


\begin{lem}
Let $P$ be a point in $\F$ such that $Du(P) = 0$. Then for any $r > 0$, there exists a connected component $U_1$ of $U$ such that $U_1 \subset B_r(P)$.
\end{lem}

\begin{proof}
First, we show that there exists a number $r' > 0$ such that if $U_1$ is any connected component of $U$ with $U_1 \cap \compl{B_r} \neq \emptyset$, then $U_1 \subset B_{r'}^c$. Indeed if it is not the case, then for any $k \in \Z$ such that $2^k < r$, there exists a connected component $U_1$ of $U$ such that
\[
U_1 \cap \compl{B_{2^k}(P)} \neq \emptyset \text{ and } U_1 \cap B_{2^{k-1}}(P) \neq \emptyset.
\]
Let $\F_1$ be a connected component of $\p U_1$ such that $\F_1$ surrounds $U_1$. We must have then that
\[
\F_1 \cap \compl{B_{2^k}(P)} \neq \emptyset \text{ and } \F_1 \cap B_{2^{k-1}}(P) \neq \emptyset.
\]
Arguing as in the Lemma \ref{t_regularity_whole_component}, we can derive the existence of a regular curve in
\[
(B_{2^k}(P) \setminus \closure{B_{2^{k-1}}(P)}) \cap \F_1
\]
and of length at least $2^k/3$. Since we can do it for all $k$ such that $2^k < r$, from the Lemma \ref{t_regularity_criterion_2} we have $\abs{Du(P)} > 0$, contradicting our hypothesis on $P$. The existence of $r'$ follows then.

Because $P \in \F$, there must exist a connected component $U_1$ of $U$ such that
\[
U_1 \cap B_{r'}(P) \neq \emptyset.
\]
The result above then guarantees that $U_1 \subset B_r(P)$.
\end{proof}


\begin{thm}
$\abs{Du} > 0$ on $\p \set{u > c}$.
\end{thm}

\begin{proof}
Assume that $\F$ contains some point where $Du = 0$. Without loss of generality, let's assume that that point is the origin. Clearly this point is not on the boundary of any connected component of $U$, as a consequence of our result in section \ref{s_regularity_component_U}.

Pick $r_1 > 0$ such that $u > c/2$ in the set $B_{r_1}$. From the previous lemma, there exists a connected component $U_1$ of $U$ such that $U_1 \subset B_{r_1}$. Since $0 \not\in \p U_1$, there exists $r_2 > 2$ such that $B_{r_2} \subset \compl{U_1}$. Choose a connected component $U_2$ of $U$ such that $U_2 \subset B_{r_2}$. Repeating for each $k$ we find a number $r_k > 0$ and a connected component $U_k$ of $U$. Let $\F_k$ be the connected component of $\p U_k$ that surrounds $U_k$. Clearly, $\F_k$ is a regular curve and the convex hull of $I(\F_k)$ is inside $B_{r_1}$. We have from definitions and the Lemma \ref{t_lower_bound_int_1_Du} that
\begin{gather}
r_1 > r_2 > \cdots \to 0\\
\closure{\F_k} \subset \F^* \cap B_{r_k} \setminus \closure{B_{r_{k+1}}}\\
\sum_{k=1}^\infty \int_{\F_k}\frac{1}{\abs{Du}} > \sum_{k=1}^\infty \frac{1}{C_1} = \infty.
\end{gather}
Applying the Lemma \ref{t_regularity_criterion} we reach a contradiction.

Thus, $\abs{Du} > 0$ on $\p U$.
\end{proof}


We combine all our results into the following statement.
\begin{thm}
Let $\Om \subset \R^2$ with Lipschitz boundary, $0 < A < \abs{\Om}$ and $\al < \closure{\al}$. Let $(u,D)$ be a minimizing configuration. Then the set $\set{u > c}$ consists of a finite number of connected components whose closures are disjoint. The boundary of each of these connected components consists of finitely many disjoint closed and simple real-analytic curves on which $\abs{Du} > 0$. Moreover, $u$ is analytic in $\closure{U}$. We can also construct a set $\tilde{D}$ such that $\p \tilde{D} = \p U$ and $\tilde{D}$, $D$ differ only in a zero measure set.
\end{thm}

\begin{proof}
Assume that there is an infinite number of connected components of $U$. Choose a sequence of distinct connected components $U_i$ of $U$ and let $P_i$ be a maximum point of $u$ in $U_i$. Let $P$ be an accumulating point of $\set{P_i}$. Because $Du(P_i) = 0$ and $u \in C^{1,1}(\Om)$, we have $Du(P) = 0$. It is trivial that $u(P) \geq c$. Now if $u(P) > c$, it means that $P$ belongs to some connected components of $U$, contradicting the fact that each $P_i$ belongs to a different connected component. So $P \in \F$ and $Du(P) = 0$, contradicting our last lemma.

Let $P$ be any point on $\p U$. Because $\abs{Du(P)} > 0$, there exists some $r > 0$ such that the set $B_r(P) \cap \set{u > c}$ is connected. Hence, $P$ is the boundary point of one and only one connected component of $U$. In other words, the closures of any two connected components do not intersect.

Assume $U_1$ is a connected component of $U$ such that $\p U_1$ consists of infinitely many connected components. Choose a sequence $\set{P_k}$ such that each $P_k$ belongs to a connected components $\F_k$ of $\p U_1$ and all $\F_k$ are distinct. Let $P \in \p U_1$ be a limit point of $\set{P_k}$. Because $\abs{Du(P)} > 0$, there exists $r > 0$ such that $B_r(P) \cap \F$ is a simple analytic curve and so, it must belongs to some connected component of $\p U_1$, contradicting the fact that each $P_k$ belongs to a different component. Thus the boundary of each connected component of $U$ consists of only a finite number of connected components.

The fact that $u$ is analytic in $\closure{U}$ is clear since $\p U$ is real-analytic, $u = c$ on $\p U$ and in $U$, $u$ satisfies the equation
\[
-\lap u = \la u.
\]

For the existence of $\tilde{D}$, just define $\tilde{D} = \compl{\closure{U}}$ and note that $\abs{\set{u = c}} = 0$.
\end{proof}

\def\cprime{$'$}
\providecommand{\bysame}{\leavevmode\hbox to3em{\hrulefill}\thinspace}
\providecommand{\MR}{\relax\ifhmode\unskip\space\fi MR }
\providecommand{\MRhref}[2]{%
  \href{http://www.ams.org/mathscinet-getitem?mr=#1}{#2}
}
\providecommand{\href}[2]{#2}

\end{document}